

\documentstyle{amsppt}
\define\z{{\Bbb Z}}
\define\cc{{\Bbb C}}
\define\r{{\Bbb R}}

\define\rk{\text{rk}\ }

\define\id{\text{id}}
\define\Tor{\text{Tor}\ }

\magnification=\magstep1

\topmatter

\title
On The Topological Classification Of Real Enriques Surfaces. I
\endtitle

\author
Viacheslav V. Nikulin
\endauthor

\dedicatory
To Memory of Dimitriy Andreevich Gudkov
\enddedicatory

\address
Steklov Mathematical Institute,
ul.Vavilova 42, Moscow 117966, GSP-1, Russia.
\endaddress

\email
slava\@nikulin.mian.su
\endemail

\abstract
This note contains our preliminary calculation of topological types
or real Enriques surfaces.
We realize 59 topological types of real Enriques surfaces (Theorem 6)
and show that all other topological types belong to the list of
21 topological types (Theorem 7). In fact, our calculation contains
much more information which is probably useful to constract or
prohibit unknown topological types.
\endabstract

\endtopmatter

\document

These notes contain our preliminary results on the topological
classification of real Enriques surfaces.

These calculations were
made during author's stay at Bielefeld University,
May---June 1992. The author is grateful to this University and
especially to Professor Heinz Helling for hospitality.
These notes were prepaired during author's stay at Kyoto
University, November 1992---January 1993. The author is
grateful to Kyoto University and especially to Professor
Masaki Maruyama for hospitalilty.

The author is grateful
to Professor R. Sujatha for useful discussions.

\smallpagebreak

The topological classification of real Enriques surfaces is a
part of the general problem of topological classification of real
algebraic varieties. See survey of D. A. Gudkov \cite{Gu}
about this problem.

\smallpagebreak

We use the following notation for the topological type of
a compact surface. We denote by $T_g$ an orientable surface of
genus $g$, and by $U_g$ a non-orientable surface of genus $g$.
Thus, $\chi (T_g)=2-2g$ and $\chi (U_g)=1-g$. The 2-sheeted
unramified orientizing covering of the $U_g$ is the surface
$T_g$.

Let $Z$ be a smooth projective algebraic surface over $\r$.
Then $Z(\cc)$ is a compact complex manifold of the complex
dimension 2, and $Z(\r)$ is a compact surface with connected
components of the type above. If $q$ is the antiholomorphic
involution on $X(\cc)$ defined by $G=\text{Gal}(\cc/\r)=\{ 1,q\}$, then
$X(\r)=X(\cc)^q$ is the set of points fixed by $q$.

By classification of algebraic surfaces (for example, see \cite{A}),
algebraic surfaces over $\cc$ are devided on the following types:
(1) General type; (2) Elliptic surfaces; (3) Abelian surfaces;
(4) Hyperelliptic surfaces; (5) K3-surfaces; (6) Enriques surfaces;
(7) Ruled surfaces; (8) Rational surfaces.

 We say that a real
surface $Z/\r$ has a type above if $Z\otimes \cc$ has this type.
Thus, we can speak about the topological classification of real
surfaces of the type above. First two classes (1) and (2) are too difficult
for classification. But classes (3)---(8) are not so large, and
their topological classification is now known (for example, see
\cite{Si}), except the type (6) of real Enriques surfaces. In the
book of R. Silhol \cite{Si}, only the topological type
$U_{10}\amalg U_0$ of real Enriques surfaces
was constructed, and we don't know some other publications on
this subject. We apply \cite{N3} and \cite{N4}
to get some new results for this classification. In these papers
some general methods of working with real K3-surfaces with
a condition on Picard lattice were developped. We also apply
results of \cite{N-S} and \cite{N5}, where some general results
on real Enriques surfaces were obtained.

First, we shortly describe our method of classification.
At the end of this notes, we give results of this classification.
Unfortunately, calculations are very long and delicate, and
we hope that they are completely correct. We hope to publish detailes of
these calculations somewhere later.

By definition (for example, see \cite{C--D}), an
Enriques surface is a quotient surface
$Y_\cc=X/\{\id,\tau \}$
where $X$ is a complex
$K3$-surface and $\tau $ an algebraic involution
of $X$ without fixed points. Let $\pi :X \to
Y_\cc $ be the quotient morphism. Since $X(\cc )$
is simply-connected, $\pi :X(\cc )\to Y(\cc )$ is
the 2-sheeted universal covering with the
holomorphic involution $\tau $ of the covering.
Thus, $\pi \tau =\pi $.

Let $\theta $ be the antiholomorphic involution
of the complex surface $Y(\cc )$ corresponding to
the real surface $Y$.
Below we use some results from \cite{N--S}.

Since $X(\cc)$ is
simply-connected, one can easily see that
there are precisely two liftings $\sigma $ and
$\tau \sigma $ of $\theta $ to
antiholomorphic automorphisms of $X(\cc )$.
If $Y(\r)\not=\emptyset$, both these
automorphisms are antiholomorphic involutions of $X(\cc)$.
Further, we suppose that this is true to consider the
empty case for $Y(\r)$ too (thus, for $Y(\r)=\emptyset$, we
suppose that these lifitings are involutions).
In other words, we lift the group
$G=\{\id,\theta\}$ of order 2 on $Y(\cc)$ to the group $\Gamma
=\{\id ,\tau ,\sigma ,\tau\sigma\}$ of order $4$
on $X(\cc)$, with $\Gamma \cong (\z/2)^2$.

We denote by $X_\sigma$ and $X_{\tau\sigma}$ the real
K3-surfaces defined by the antiholomorphic involutions
$\sigma$ and $\tau\sigma$ respectively. Thus,
$X_\sigma (\r)=X(\cc)^\sigma$ and
$X_{\tau\sigma}(\r)=X(\cc)^{\tau\sigma}$.
Becides, $\tau (X_\sigma(\r))=X_\sigma(\r)$,
$\tau (X_{\tau\sigma} (\r))=X_{\tau\sigma} (\r)$.
The $Y(\r)=\pi (X_\sigma (\r))\amalg \pi (X_{\tau\sigma }(\r))$,
and we denote
$Y(\r)_\sigma =\pi (X_\sigma (\r))=X_\sigma (\r)/\{ 1, \tau\}$,
$Y(\r)_{\tau\sigma}=$
\linebreak
$\pi (X_{\tau\sigma}(\r))=X_{\tau\sigma}(\r)/\{ 1, \tau\}$.

Let $\omega _X$ be a non-zero holomorphic 2-form of $X$.
The corresponding real parts
$\omega _X^\sigma$ and $\omega _X^{\tau\sigma}$ of
$\omega _X$ define the  canonical volume forms on
$X_\sigma (\r)$ and $X_{\tau\sigma} (\r)$. Since
$\tau (\omega _X)=-\omega _X$ (it is well-known),
the canonical map
$$
\pi : X_\sigma (\r) \amalg X_{\tau\sigma} (\r) \to
Y(\r)_\sigma \amalg Y(\r)_{\tau\sigma}=Y(\r)
\tag1
$$
is the 2-sheeted unramified orientizing covering.

Thus, we can correspond to a real Enriques surface $Y$ the
topological type of the five surfaces
$$
(X_\sigma (\r), X_{\tau\sigma}(\r), Y(\r)_\sigma,
Y(\r)_{\tau\sigma}, Y(\r)=Y(\r)_\sigma \amalg
Y(\r)_{\tau\sigma}).
\tag2
$$
We want to get an information about these five surfaces.

Here $X_\sigma $ and $X_{\tau\sigma}$ are real K3-surfaces. Thus,
at first, we recall the theory of real K3-surfaces.

Let $X$ be a real K3-surfaces with the antiholomorphic involution
$\phi $. We recall that the cohomology lattice (with the intersection
pairing) $H^2(X(\cc);\z)$ is isomorphic to
a standard even unimodular lattice $L$
of signature $(3,19)$. Here we have the following result which
follows from Global Torelli Theorem \cite{P\u S-\u S} and
Epimorphicity of the Torelli map \cite{Ku} for K3-surfaces,
and the geometrical interpretation of invariants.

\proclaim{Theorem 1 (Kharlamov \cite{Ha1}, Nikulin \cite{N3})} The
topological type $X(\r)=X(\cc )^\phi$ is defined by the
action of $\phi$ on the cohomology lattice $H^2(X(\cc);\z)$.

An action of the abstract group $G=\{1, \phi \}$ of order two
on the lattice $L$ corresponds to a real K3-surface if and only if
the lattice $L^\phi$ is hyperbolic (it has the signature $(1, t)$).
\endproclaim

All actions of $G=\{1, \phi \}$ are classified by some invariants
$$
(r(\phi ), a(\phi ), \delta (\phi))
\tag3
$$
(see \cite{N3}):

Let $L_+=L^\phi $ and $L_-=L^{(-\phi )}$. Then we have
an orthogonal decomposition up to a finite index
$L_+ + L_-\subset L$. The invariant $r(\phi )=\rk L_+$.

We have
$$L_+/2L_+\subset L/2L \supset L_-/2L_-$$
$$\cup \ \ \ \ \ \ \ \ \ \cup \ \ \ \ \ \ \ \ \ \cup
\tag4$$
$$A(\phi )_+ = A(\phi )=A(\phi )_-$$
where $A(\phi )=L_+/2L_+ \cap L_-/2L_-$. Here
$A(\phi )=(\z/2\z)^{a(\phi)}$. This gives the invariant
$a(\phi)$.

Let
$$
q(\phi)(v_+ + 2L_+)=(1/2)v_+^{\ 2}\mod 4,
\tag5
$$
for $v_+ + 2L_+ \in A(\phi )_+$.
This defines a quadratic form $q(\phi)$ on $A(\phi)$ with the
non-degenerate symmetric bilinear form $b(\phi)$:
$$
b(\phi)(v_+ + 2L_+, w_+ + 2L_+)=(1/2)v_+\cdot w_+\mod 2.
\tag6
$$
The invariant $\delta (\phi)=0$ if
$q(\phi)(A(\phi))\subset 2\z /4\z$. Otherwise, $\delta (\phi)=1$.

There exists another definition of the invariant $\delta (\phi)$.
The element $s_\phi \in L/2L$ is called characteristic if
$$
x\cdot \phi (x)\equiv s_\phi \cdot x \mod 2
\tag7
$$
for any $x\in L$. One can see that $\delta (\phi )=0$ iff the
characteristic element $s_\phi$ is equal to $0$.

Using these invariants, we have
$$
X(\cc )^\phi =
\cases
\emptyset, & \text{if $(r(\phi),a(\phi),\delta(\phi))=(10,10,0)$};\\
2T_1, & \text{if $(r(\phi), a(\phi),\delta(\phi))=(10,8,0)$};\\
T_{g(\phi )}\coprod k(\phi)T_0 &
\text{where $g(\phi)=(22-r(\phi)-a(\phi))/2$},\\
& \text{$k(\phi)=(r(\phi)-a(\phi))/2$ otherwise}.
\endcases
\tag8
$$
We remark that
$$
\chi(X(\r ))=2r(\phi )-10;
\tag9
$$
$$
\dim H_\ast (X(\r );\z/2)=24-2a(\phi )\ \text{if} \ X(\r)\not=\emptyset .
\tag10
$$
The characteristic element
$$
s_\phi \equiv  X(\r ) \mod 2
$$
in $H_2(X(\cc);\z)$ (here, we identify cohomology and homology using
Poincar\'e duality). Thus,
$$
X(\r )\equiv 0 \mod 2\ \text{iff} \ \delta (\phi )=0.
\tag12
$$
All possibilities for triplets (3) are known (see \cite{N3}), and this
gives the topological classification of real $K3$-surfaces (see
\cite{Ha1} and \cite{N3}).

We want to apply similar idea to study real Enriques surfaces, but for
the action of the group
$$
\Gamma =\{ \text{id}, \tau, \sigma, \tau\sigma\} \cong (\z /2)^2
\tag13
$$
on $L$.

Here the action of $\tau$ is standard. It is defined uniquely by the
condition that
$$
L^\tau =S \cong E(2)
$$
where $E$ is an even unimodular lattice of signature $(1,9)$, and
$E(2)$ means that we multiply the form of $E$ by $2$.
The lattice $E=H^2(Y(\cc ); \z )/ \Tor $, and $S=\pi ^\ast E$.
Thus, the action
$$
\sigma \mid S=\tau\sigma \mid S \cong \theta \mid E.
$$

By global Torelli Theorem \cite{P\u S-\u S} and epimorphisity of the
Torelli map \cite{Ku}, we have similar to Theorem 1

\proclaim{Theorem 2} An action of the group
$\Gamma =\{\text{id}, \tau, \sigma, \tau\sigma \} \cong (\z/2)^2$ on
$L$ corresponds to a real Enriques surface iff
$L^\tau= S\cong E(2)$ and the lattices $L^\sigma$ and
$L^{\tau\sigma}$ are both hyperbolic.
\endproclaim

Unfortunately, we don't know that the action of $\Gamma $ defines the
topology of a five (2). We only have particular results about.

Thus, our problem of topological classification of real Enriques
surfaces is devided on two:

(A) To classify actions of $\Gamma $ on $L$
with conditions of Theorem 2.

(B) To find the geometrical interpretation of invariants of the actions of
$\Gamma $ on $L$.

The problem (A):

First, we should classify actions $\sigma \mid S=\tau\sigma \mid S$
which are equivalent to the action $\theta \mid E$. Here the lattice
$E$ is unimodular, and we have similar to (3) invariants
$$
(r(\theta), a(\theta), \delta (\theta))
\tag15
$$
for this action. Here the lattice $E^\theta $ is negative definite.
Using results of \cite{N3} or \cite{N4}, we have the following
possibilities for these invariants:
$$
\split
(r(\theta ), a(\theta ), \delta (\theta ))=
&(1,1,1),\ (2,2,1),\ (3,3,1),\ (4,4,1),\ (5,5,1), \\
& (9,1,1),\ (8,2,1),\ (7,3,1),\ (6,4,1),\ (0,0,0),\ (8,2,0),\\
&(5,3,1),\ (6,2,1),\ (7,1,1),\ (4,2,0),\ (8,0,0).
\endsplit
\tag16
$$
We fix one of these possibilities.

Second, we choose
an involution $\sigma$ between $\sigma $ and $\tau\sigma $.
In fact, one should choose between invariants
$(r(\sigma), a(\sigma), \delta (\sigma))$ and
$(r(\tau\sigma ), a(\tau\sigma), \delta (\tau\sigma ))$ of these
involutions, and almost in all cases this choice is canonical.
The formulae (24), (25) and (26) below are useful for.

Third, we should study the problem of extension of
$\sigma \mid S$ to  the unimodular lattice $L$. The paper \cite{N4}
was devoted to this problem. In \cite{N4}, there were found
all invariants of genus
of these extensions for an arbitrary lattice $S$. Becides, all
relations between these invariants which are necessary and
sufficient for existence of these extensions were found too. We need
to describe these invariants here.

We have the canonical subgroups $H(\sigma )_+, H(\sigma )_-$
(see (4)):
$$
(S_+\cap S_-)/2(S_+\cap S_-)\subset
H(\sigma )_\pm=S_\pm/2S_\pm \cap A(\sigma )_\pm \subset S_\pm /2S_\pm .
\tag17
$$
Using identifications $A(\sigma )_\pm =A(\sigma )$, we can
consider these subgroups $H(\sigma)_\pm$ as subgroups in
the space $A(\sigma)$ equipped with the bilinear form $b(\sigma )$.
This defines the canonical pairing
$$
\rho (\sigma ): H(\sigma )_+ \times H(\sigma )_- \to \z/2\z.
\tag18
$$
The characteristic vector $s_\sigma$ may belong to $S/2S$,
$(S_+\cap S_-)/2(S_+\cap S_-)$ or not.
This defines similar to $\delta (\sigma )$
invariants $\delta_{\sigma S}$, $\delta_ {\sigma S_+\cap S_-}$,
which are equal to $0$ or $1$. For example,
$\delta_{\sigma S}=0$ iff the $s_\sigma \in S/2S$.
If $\delta_{\sigma S}=0$, the element
$s_\sigma \in S/2S$ gives an additional invariant of the extension.

It was proved in \cite{N4}, that the invariants
$$
r(\sigma ), a(\sigma ), \delta (\sigma ),
H(\sigma )_+, H(\sigma )_-, \delta _{\sigma S},
s_\sigma \in S/2S\  (if\ \delta _{\sigma S}=0)
\tag19
$$
together with real invariants of $\sigma$ on $L$ are the all invariants
of the genus of extensions of $\sigma$ on
even unimodular lattices $L$.
Using necessary and sufficient conditions from \cite{N4}, we can
describe all invariants (19) for real Enriques surfaces.

Here,  we will restict by the most important for topology invariants
\newline
$(r(\sigma ), a(\sigma ), \delta (\sigma ))$ and
$$
h(\sigma )_\pm =\dim H(\sigma )_\pm,
\tag20
$$
$$
\gamma (\sigma )=\dim H(\sigma )_- - \dim (H(\sigma )_+)^\perp \cap
H(\sigma )_-,
\tag21
$$
(here we use the pairing $\rho (\sigma )$), and
$$
\delta_{\sigma S}, \delta_{\sigma S_+\cap S_-}.
\tag22
$$

Using \cite{N4} (see \cite{N--S}), we can prove that
$$
a(\sigma )=h(\sigma )_+ + h(\sigma )_- +\alpha (\sigma )
\tag23
$$
where the invariant $\alpha (\sigma)=0$ or $1$.

\smallpagebreak

The Problem (B).

First of all, we have the following formulae, which connect the
invariants above with the invariants
$(r(\tau\sigma), a(\tau\sigma ), \delta (\tau\sigma))$
(see \cite{N--S}).
$$
r(\sigma)+r(\tau\sigma)=12+2r(\theta);
\tag24
$$
$$
a(\sigma)+a(\tau\sigma)=10+2a(\theta)+2\gamma(\sigma)+2\alpha (\sigma);
\tag25
$$
$$
\delta (\sigma)+\delta (\tau\sigma)\equiv \delta (\theta) \mod 2
\tag26
$$
(recall (12)).
Thus, using the invariants (19), we can find the invariants
\linebreak
$(r(\tau\sigma), a(\tau\sigma ), \delta (\tau\sigma))$
of $\tau\sigma$. Also, by formulae above, the invariants
\linebreak
$(r(\theta), a(\theta), \delta (\theta))$, and
$\alpha (\sigma )$, $\delta (\sigma )$, $\delta (\tau\sigma )$
and the topological type of $X_\sigma (\r )$ and $X_{\tau\sigma } (\r )$
define the invariants
(3) for $\sigma $ and $\tau\sigma$, and (20), (21). Thus, to formulate
results of calculations, we can use these invariants instead of (20) and
(21).

The very non-trivial formula is a formula for the number
$s_{nor}$ of non-orientable connected components of $Y(\r)$ for
a real Enriques surface $Y$. These uses results of
\cite{N--S} and
\cite{N5}.

\proclaim{Theorem 3} Let both $X_\sigma (\r)$ and
$X_{\tau\sigma}(\r)$ are non-empty.

Then
$$
s_{nor}=1 +
\alpha(\sigma )(2\delta_{\sigma S_+\cap S_-} - 1)+
\gamma (\sigma) .
\tag3--6
$$
\endproclaim

When $X_\sigma (\r)=\emptyset$, we only have an inequality

\proclaim{Theorem (see \cite{N--S}) 4} Let $X_\sigma (\r)=\emptyset$
(it follows that $\delta_{\sigma S_+\cap S_-}=0$)
and $X_{\tau\sigma }(\r)\not=\emptyset$.

Then
$$
s_{nor}\le 2-\alpha (\sigma )+\gamma (\sigma )
$$
\endproclaim

To use these Theorems, we add the invariant
$\delta_{\sigma S_+\cap S_-}$ to invarants above. Evidently,
$\delta _{\sigma S}=\delta_{\sigma S_+\cap S_-}=0$ if
$\delta (\sigma )=0$, and
$\delta (\sigma ) = \delta_{\sigma S_+\cap S_-}=1$ if $\delta _{\sigma S}=1$.
In many cases, Theorems 3 and 4 permit to find the topological type of
$Y(\r)$ and $Y(\r)_\sigma, Y(\r)_{\tau\sigma}$ using the orientizing
map (1).

Below, we give results of calculations using the method described above,
where "A or B" means
that we are proving that
one of possibilities A, B, (A and B) holds, but we don't know
what of them does hold.

\proclaim{Theorem 5}
We have the following and only following possibilities for the
\linebreak
invariants
\newline
$(r(\theta), a(\theta), \delta(\theta));\
\alpha (\sigma );\ X_\sigma (\r),
\delta (\sigma ), \delta_{\sigma S}, \delta_{\sigma S_+\cap S_-};
 X_{\tau\sigma} (\r),
\delta (\tau\sigma ), \delta_{\tau\sigma S}, \delta_{\tau\sigma S_+\cap S_-}$;
\linebreak
$Y(\r)_\sigma , Y(\r)_{\tau\sigma}, Y(\r)$ for real Enriques surfaces $Y$:
\endproclaim

\bigpagebreak

The case $(r(\theta), a(\theta), \delta (\theta))=(1,1,1)$:

\smallpagebreak

\centerline
{$
\alpha (\sigma)=1,
X_\sigma (\r)=\emptyset,
\delta (\sigma)=0,
X_{\tau\sigma}(\r)=T_7,
\delta (\tau\sigma )=1,
$}
\centerline{
$
Y(\r)_\sigma=\emptyset,
Y(\r)_{\tau\sigma}=U_7,
Y(\r)=U_7;
$}

\smallpagebreak

\centerline{
$
\alpha (\sigma)=0,
X_\sigma (\r)=\emptyset,
\delta (\sigma)=0,
X_{\tau\sigma}(\r)=T_8\amalg T_0,
\delta (\tau\sigma )=1,
$}
\centerline{
$
Y(\r)_\sigma=\emptyset,
Y(\r)_{\tau\sigma}=U_8\amalg U_0,
Y(\r)=U_8\amalg U_0;
$}

\smallpagebreak

\centerline{
$
\alpha (\sigma)=0,
X_\sigma (\r)=2T_1,
\delta (\sigma)=0,
X_{\tau\sigma}(\r)=T_7,
\delta (\tau\sigma )=1,
$}
\centerline{
$
Y(\r)_\sigma=T_1,
Y(\r)_{\tau\sigma}=U_7,
Y(\r)=U_7\amalg T_1;
$}

\smallpagebreak

\centerline{
$
\alpha (\sigma)=0,
X_\sigma (\r)=T_9,
\delta (\sigma)=0,
X_{\tau\sigma}(\r)=2T_0,
\delta (\tau\sigma )=1,
$}
\centerline{
$
Y(\r)_\sigma=U_9,
Y(\r)_{\tau\sigma}=T_0,
Y(\r)=U_9\amalg T_0.
$}

\bigpagebreak

The case $(r(\theta), a(\theta ), \delta (\theta ))=(2,2,1)$:

\smallpagebreak

\centerline{
$
\alpha (\sigma)=1,
X_\sigma (\r)=\emptyset,
\delta (\sigma)=0,
X_{\tau\sigma}(\r)=T_5,
\delta (\tau\sigma )=1,
$}
\centerline{
$
Y(\r)_\sigma=\emptyset,
Y(\r)_{\tau\sigma}=U_5,
Y(\r)=U_5;
$}

\smallpagebreak

\centerline{
$
\alpha (\sigma)=0,
X_\sigma (\r)=\emptyset,
\delta (\sigma)=0,
X_{\tau\sigma}(\r)=T_6\amalg T_0,
\delta (\tau\sigma )=1,
$}
\centerline{
$
Y(\r)_\sigma=\emptyset,
Y(\r)_{\tau\sigma}=U_6\amalg U_0,
Y(\r)=U_6\amalg U_0;
$}

\smallpagebreak

\centerline{
$
\alpha (\sigma)=0;
X_\sigma (\r)=2T_1,
\delta (\sigma)=0,
X_{\tau\sigma}(\r)=T_5,
\delta (\tau\sigma )=1,
$}
\centerline{
$
Y(\r)_\sigma=T_1,
Y(\r)_{\tau\sigma}=U_5,
Y(\r)=U_5\amalg T_1.
$}

\bigpagebreak

The case $(r(\theta), a(\theta ), \delta (\theta ))=(3,3,1)$:

\smallpagebreak

\centerline{
$
\alpha (\sigma)=1;
X_\sigma (\r)=\emptyset,
\delta (\sigma)=0,
X_{\tau\sigma}(\r)=T_3,
\delta (\tau\sigma )=1,
$}
\centerline{
$
Y(\r)_\sigma=\emptyset,
Y(\r)_{\tau\sigma}=U_3,
Y(\r)=U_3;
$}

\smallpagebreak

\centerline{
$
\alpha (\sigma)=0;
X_\sigma (\r)=\emptyset,
\delta (\sigma)=0,
X_{\tau\sigma}(\r)=T_4\amalg T_0,
\delta (\tau\sigma )=1,
$}
\centerline{
$
Y(\r)_\sigma=\emptyset,
Y(\r)_{\tau\sigma}=U_4\amalg U_0,
Y(\r)=U_4\amalg U_0;
$}

\smallpagebreak

\centerline{
$
\alpha (\sigma)=0;
X_\sigma (\r)=2T_1,
\delta (\sigma)=0,
X_{\tau\sigma}(\r)=T_3,
\delta (\tau\sigma )=1,
$}
\centerline{
$
Y(\r)_\sigma=T_1,
Y(\r)_{\tau\sigma}=U_3,
Y(\r)=U_3\amalg T_1.
$}

\bigpagebreak

The case $(r(\theta), a(\theta ), \delta (\theta ))=(4,4,1)$:

\smallpagebreak

\centerline{
$
\alpha (\sigma)=1;
X_\sigma (\r)=\emptyset,
\delta (\sigma)=0,
X_{\tau\sigma}(\r)=T_1,
\delta (\tau\sigma )=1,
$}
\centerline{
$
Y(\r)_\sigma=\emptyset,
Y(\r)_{\tau\sigma}=U_1,
Y(\r)=U_1;
$}

\smallpagebreak

\centerline{
$
\alpha (\sigma)=0;
X_\sigma (\r)=\emptyset,
\delta (\sigma)=0,
X_{\tau\sigma}(\r)=T_2\amalg T_0,
\delta (\tau\sigma )=1,
$}
\centerline{
$
Y(\r)_\sigma=\emptyset,
Y(\r)_{\tau\sigma}=U_2\amalg U_0,
Y(\r)=U_2\amalg U_0;
$}

\smallpagebreak

\centerline{
$
\alpha (\sigma)=0;
X_\sigma (\r)=2T_1,
\delta (\sigma)=0,
X_{\tau\sigma}(\r)=T_1,
\delta (\tau\sigma )=1,
$}
\centerline{
$
Y(\r)_\sigma=T_1,
Y(\r)_{\tau\sigma}=U_1,
Y(\r)=U_1\amalg T_1.
$}

\bigpagebreak

The case $(r(\theta), a(\theta ), \delta (\theta ))=(5,5,1)$:

\smallpagebreak

\centerline{
$
\alpha (\sigma)=0;
X_\sigma (\r)=\emptyset,
\delta (\sigma)=0,
X_{\tau\sigma}(\r)=2T_0,
\delta (\tau\sigma )=1,
$}
\centerline{
$
Y(\r)_\sigma=\emptyset,
Y(\r)_{\tau\sigma}=Y(\r)= 2U_0\ or\ T_0.
$}

\bigpagebreak

The case $(r(\theta), a(\theta ), \delta (\theta ))=(9,1,1)$:

\smallpagebreak

\centerline{
$
\alpha (\sigma)=0;
X_\sigma (\r)=\emptyset,
\delta (\sigma)=0,
X_{\tau\sigma}(\r)=10T_0,
\delta (\tau\sigma )=1,
$}
\centerline{
$
Y(\r)_\sigma=\emptyset,
Y(\r)_{\tau\sigma}=Y(\r)=5T_0\ or\ 2U_0\amalg 4T_0;
$}

\smallpagebreak

\centerline{
$
\alpha (\sigma)=1;
X_\sigma (\r)=8T_0,
\delta (\sigma)=0,
X_{\tau\sigma}(\r)=2T_0,
\delta (\tau\sigma )=1,
$}
\centerline{
$
Y(\r)_\sigma=4T_0,
Y(\r)_{\tau\sigma}=T_0,
Y(\r)=5T_0;
$}

\smallpagebreak

\centerline{
$
\alpha (\sigma)=0;
X_\sigma (\r)=T_1\amalg 8T_0,
\delta (\sigma)=0,
X_{\tau\sigma}(\r)=2T_0,
\delta (\tau\sigma )=1,
$}
\centerline{
$
Y(\r)_\sigma=U_1\amalg 4T_0,
Y(\r)_{\tau\sigma}=T_0,
Y(\r)=U_1\amalg 5T_0;
$}

\smallpagebreak

\centerline{
$
\alpha (\sigma)=0;
X_\sigma (\r)=8T_0,
\delta (\sigma)=0,
X_{\tau\sigma}(\r)=T_1\amalg 2T_0,
\delta (\tau\sigma )=1,
$}
\centerline{
$
Y(\r)_\sigma=4T_0,
Y(\r)_{\tau\sigma}=U_1\amalg T_0,
Y(\r)=U_1\amalg 5T_0.
$}

\bigpagebreak

The case $(r(\theta), a(\theta ), \delta (\theta ))=(8,2,1)$:

\smallpagebreak

\centerline{
$
\alpha (\sigma)=1;
X_\sigma (\r)=4T_0,
\delta (\sigma)=0,
X_{\tau\sigma}(\r)=4T_0,
\delta (\tau\sigma )=1,
$}
\centerline{
$
Y(\r)_\sigma=2T_0,
Y(\r)_{\tau\sigma}=2T_0,
Y(\r)=4T_0;
$}

\smallpagebreak

\centerline{
$
\alpha (\sigma)=0;
X_\sigma (\r)=T_1\amalg 4T_0,
\delta (\sigma)=0,
X_{\tau\sigma}(\r)=4T_0,
\delta (\tau\sigma )=1,
$}
\centerline{
$
Y(\r)_\sigma=U_1\amalg 2T_0,
Y(\r)_{\tau\sigma}=2T_0,
Y(\r)=U_1\amalg 4T_0;
$}

\smallpagebreak

\centerline{
$
\alpha (\sigma)=0;
X_\sigma (\r)=4T_0,
\delta (\sigma)=0,
X_{\tau\sigma}(\r)=T_1\amalg 4T_0,
\delta (\tau\sigma )=1,
$}
\centerline{
$
Y(\r)_\sigma=2T_0,
Y(\r)_{\tau\sigma}=U_1\amalg 2T_0,
Y(\r)=U_1\amalg 4T_0.
$}

\bigpagebreak

The case $(r(\theta), a(\theta ), \delta (\theta ))=(7,3,1)$:

\smallpagebreak

\centerline{
$
\alpha (\sigma)=0;
X_\sigma (\r)=\emptyset,
\delta (\sigma)=0,
X_{\tau\sigma}(\r)=6T_0,
\delta (\tau\sigma )=1,
$}
\centerline{
$
Y(\r)_\sigma=\emptyset,
Y(\r)_{\tau\sigma}=Y(\r)=3T_0\ or\ 2U_0\amalg 2T_0.
$}

\bigpagebreak

The case $(r(\theta), a(\theta ), \delta (\theta ))=(6,4,1)$:

\smallpagebreak

\centerline{
$
\alpha (\sigma)=0;
X_\sigma (\r)=4T_0,
\delta (\sigma)=0,
X_{\tau\sigma}(\r)=T_1,
\delta (\tau\sigma )=1,
$}
\centerline{
$
Y(\r)_\sigma =2T_0,
Y(\r)_{\tau\sigma}=U_1,
Y(\r)=U_1\amalg 2T_0.
$}

\bigpagebreak

The case $(r(\theta), a(\theta ), \delta (\theta ))=(0,0,0)$:

\smallpagebreak

\centerline{
$
\alpha (\sigma)=1;
X_\sigma (\r)=T_9,
\delta (\sigma)=0,
X_{\tau\sigma}(\r)=\emptyset,
\delta (\tau\sigma )=0,
$}
\centerline{
$
Y(\r)_\sigma=U_9,
Y(\r)_{\tau\sigma}=\emptyset,
Y(\r)=U_9;
$}

\smallpagebreak

\centerline{
$
\alpha (\sigma)=0;
X_\sigma (\r)=T_{10}\amalg T_0,
\delta (\sigma)=0,
X_{\tau\sigma}(\r)=\emptyset,
\delta (\tau\sigma )=0,
$}
\centerline{
$
Y(\r)_\sigma=U_{10}\amalg U_0,
Y(\r)_{\tau\sigma}=\emptyset,
Y(\r)=U_{10}\amalg U_0;
$}

\smallpagebreak

\centerline{
$
\alpha (\sigma)=0;
X_\sigma (\r)=T_9,
\delta (\sigma)=0,
X_{\tau\sigma}(\r)=2T_1,
\delta (\tau\sigma )=0,
$}
\centerline{
$
Y(\r)_\sigma=U_9,
Y(\r)_{\tau\sigma}=T_1,
Y(\r)=U_9\amalg T_1;
$}

\smallpagebreak

\centerline{
$
\alpha (\sigma)=1;
X_\sigma (\r)=T_{10-2t},
\delta_{\sigma S}=1,
X_{\tau\sigma}(\r)=T_{2t},
\delta _{\tau\sigma S}=1,
$}
\centerline{
$
Y(\r)_\sigma =U_{10-2t},
Y(\r)_{\tau\sigma}=U_{2t},
Y(\r)=U_{10-2t}\amalg U_{2t},\ where\ t=0,1,2;
$}

\smallpagebreak

\centerline{
$
\alpha (\sigma)=1;
X_\sigma (\r)=T_{9-2t},
\delta_{\sigma S}=0,
\delta_{\sigma S_+ \cap S_-} = 1,
$}
\centerline{
$X_{\tau\sigma}(\r)=T_{1+2t},
\delta _{\tau\sigma S}=0,
\delta_{\tau \sigma S_+ \cap S_-} = 1,
$}
\centerline{
$
Y(\r)_\sigma =U_{9-2t},
Y(\r)_{\tau\sigma}=U_{1+2t},
Y(\r)=U_{9-2t}\amalg U_{1+2t},\ where\ t=0,1,2.
$}

\bigpagebreak

The case $(r(\theta), a(\theta ), \delta (\theta ))=(8,2,0)$:

\smallpagebreak

\centerline{
$
\alpha (\sigma)=0;
X_\sigma (\r)=\emptyset,
\delta (\sigma)=0,
X_{\tau\sigma}(\r)=8T_0,
\delta (\tau\sigma )=0,
$}
\centerline{
$
Y(\r)_\sigma=\emptyset,
Y(\r)_{\tau\sigma}=Y(\r)=4T_0\ or\ 2U_0\amalg 3T_0;
$}

\smallpagebreak

\centerline{
$
\alpha (\sigma)=1;
X_\sigma (\r)=(1+2t)T_0,
\delta _{\sigma S}=1,
X_{\tau\sigma}(\r)=(7-2t)T_0,
\delta_{\tau\sigma S} = 1,
$}
\centerline{
$
Y(\r)_\sigma=U_0\amalg tT_0,
Y(\r)_{\tau\sigma}=U_0\amalg (3-t)T_0,
Y(\r)=2U_0\amalg 3T_0,\ where\ t=0,1;
$}

\smallpagebreak

\centerline{
$
\alpha (\sigma)=1;
X_\sigma (\r)=2T_0,
\delta _{\sigma S}=0,
\delta _{\sigma S_+\cap S_-}=1,
$}
\centerline{
$
X_{\tau\sigma}(\r)=6T_0,
\delta _{\tau\sigma S}=0,
\delta _{\tau\sigma S_+\cap S_-}=1,
$}
\centerline{
$
(Y(\r)_\sigma,
Y(\r)_{\tau\sigma})=(2U_0,3T_0)\ or\ (T_0,2U_0\amalg 2T_0),
Y(\r)=2U_0\amalg 3T_0;
$}

\smallpagebreak

\centerline{
$
\alpha (\sigma)=1;
X_\sigma (\r)=4T_0,
\delta _{\sigma S}=0,
\delta _{\sigma S_+\cap S_-}=1,
$}
\centerline{
$
X_{\tau\sigma}(\r)=4T_0,
\delta _{\tau\sigma S}=0,
\delta _{\tau\sigma S_+\cap S_-}=1,
$}
\centerline{
$
Y(\r)_\sigma =2U_0\amalg T_0,
Y(\r)_{\tau\sigma}=2T_0,
Y(\r)=2U_0\amalg 3T_0;
$}

\smallpagebreak

\centerline{
$
\alpha (\sigma)=1;
X_\sigma (\r)=4T_0,
\delta (\sigma )=1,
\delta _{\sigma S_+\cap S_-}=0,
$}
\centerline{
$
X_{\tau\sigma}(\r)=4T_0,
\delta (\tau\sigma )=1,
\delta _{\tau\sigma S_+\cap S_-}=0,
$}
\centerline{
$
Y(\r)_\sigma =2T_0,
Y(\r)_{\tau\sigma}=2T_0,
Y(\r)=4T_0;
$}

\smallpagebreak

\centerline{
$
\alpha (\sigma)=0;
X_\sigma (\r)=T_1\amalg (4-4t)T_0,
\delta (\sigma )=1,
\delta _{\sigma S_+\cap S_-}=0,
$}
\centerline{
$
X_{\tau\sigma} (\r )=(4+4t)T_0, \delta (\tau\sigma )=1,
\delta_{\tau\sigma S_+\cap S_-}=0,
$}
\centerline{
$
Y(\r )_\sigma=U_1\amalg (2-2t)T_0, Y(\r )_{\tau\sigma}=(2+2t)T_0,
Y(\r )=U_1\amalg 4T_0,\ where\ t=0,1.
$}

\bigpagebreak

The case $(r(\theta), a(\theta ), \delta (\theta ))=(5,3,1)$:

\smallpagebreak

\centerline{
$
\alpha (\sigma)=1;
X_\sigma (\r)=\emptyset ,
\delta (\sigma)=0,
X_{\tau\sigma}(\r)=2T_0,
\delta (\tau\sigma )=1,
$}
\centerline{
$
Y(\r)_\sigma=\emptyset ,
Y(\r)_{\tau\sigma}=Y(\r)=2U_0\ or\ T_0;
$}

\smallpagebreak

\centerline{
$
\alpha (\sigma)=0;
X_\sigma (\r)=\emptyset ,
\delta (\sigma)=0,
X_{\tau\sigma}(\r)=2T_0,
\delta (\tau\sigma )=1,
$}
\centerline{
$
Y(\r)_\sigma=\emptyset ,
Y(\r)_{\tau\sigma}=Y(\r)=2U_0\ or\ T_0;
$}

\smallpagebreak

\centerline{
$
\alpha (\sigma)=0;
X_\sigma (\r)=\emptyset ,
\delta (\sigma)=0,
X_{\tau\sigma}(\r)=T_1\amalg 2T_0,
\delta (\tau\sigma )=1,
$}
\centerline{
$
Y(\r)_\sigma=\emptyset,
Y(\r)_{\tau\sigma}=Y(\r)=U_1\amalg 2U_0\ or\ U_1 \amalg T_0;
$}

\smallpagebreak

\centerline{
$
\alpha (\sigma)=1;
X_\sigma (\r)=\emptyset ,
\delta (\sigma)=0,
X_{\tau\sigma}(\r)=T_1\amalg 2T_0,
\delta (\tau\sigma )=1,
$}
\centerline{
$
Y(\r)_\sigma=\emptyset,
Y(\r)_{\tau\sigma}=Y(\r)=U_1\amalg T_0;
$}

\smallpagebreak

\centerline{
$
\alpha (\sigma)=0;
X_\sigma (\r)=\emptyset ,
\delta (\sigma)=0,
X_{\tau\sigma}(\r)=T_2\amalg 3T_0,
\delta (\tau\sigma )=1,
$}
\centerline{
$
Y(\r)_\sigma=\emptyset,
Y(\r)_{\tau\sigma}=Y(\r)=U_2\amalg U_0\amalg T_0;
$}

\smallpagebreak

\centerline{
$
\alpha (\sigma)=0;
X_\sigma (\r)=2T_1 ,
\delta (\sigma)=0,
X_{\tau\sigma}(\r)=2T_0,
\delta (\tau\sigma )=1,
$}
\centerline{
$
(Y(\r)_\sigma, Y(\r)_{\tau\sigma},Y(\r))=
(2U_1, T_0, 2U_1\amalg T_0)\ or\
(T_2, 2U_0, 2U_0\amalg T_1);
$}

\smallpagebreak

\centerline{
$
\alpha (\sigma)=1;
X_\sigma (\r)=2T_1,
\delta (\sigma)=0,
X_{\tau\sigma}(\r)=2T_0,
\delta (\tau\sigma )=1,
$}
\centerline{
$
Y(\r)_\sigma=T_1,
Y(\r)_{\tau\sigma}=T_0,
Y(\r)=T_1\amalg T_0;
$}

\smallpagebreak

\centerline{
$
\alpha (\sigma)=0;
X_\sigma (\r)=2T_1,
\delta (\sigma)=0,
X_{\tau\sigma}(\r)=T_1\amalg 2T_0,
\delta (\tau\sigma )=1,
$}
\centerline{
$
Y(\r)_\sigma=T_1,
Y(\r)_{\tau\sigma}=U_1\amalg T_0,
Y(\r)=U_1\amalg T_1\amalg T_0;
$}

\smallpagebreak

\centerline{
$
\alpha (\sigma)=0;
X_\sigma (\r)=T_3\amalg 2T_0 ,
\delta (\sigma)=0,
X_{\tau\sigma}(\r)=2T_0,
\delta (\tau\sigma )=1,
$}
\centerline{
$
Y(\r)_\sigma= U_3\amalg T_0,
Y(\r)_{\tau\sigma}=T_0,
Y(\r)=U_3\amalg 2T_0;
$}

\bigpagebreak

The case $(r(\theta), a(\theta ), \delta (\theta ))=(6,2,1)$:

\smallpagebreak

\centerline{
$
\alpha (\sigma)=1;
X_\sigma (\r)=\emptyset,
\delta (\sigma)=0,
X_{\tau\sigma}(\r)=4T_0,
\delta (\tau\sigma )=1,
$}
\centerline{
$
Y(\r)_\sigma=\emptyset,
Y(\r)_{\tau\sigma}=Y(\r)=2U_0\amalg T_0\ or\ 2T_0;
$}

\smallpagebreak

\centerline{
$
\alpha (\sigma)=0;
X_\sigma (\r)=\emptyset,
\delta (\sigma)=0,
X_{\tau\sigma}(\r)=4T_0,
\delta (\tau\sigma )=1,
$}
\centerline{
$
Y(\r)_\sigma=\emptyset,
Y(\r)_{\tau\sigma}=Y(\r)=2U_0\amalg T_0\ or\ 2T_0;
$}

\smallpagebreak

\centerline{
$
\alpha (\sigma)=0;
X_\sigma (\r)=\emptyset,
\delta (\sigma)=0,
X_{\tau\sigma}(\r)=T_1\amalg 4T_0,
\delta (\tau\sigma )=1,
$}
\centerline{
$
Y(\r)_\sigma=\emptyset,
Y(\r)_{\tau\sigma}=Y(\r)=U_1\amalg 2U_0\amalg T_0\ or\ U_1\amalg 2T_0;
$}

\smallpagebreak

\centerline{
$
\alpha (\sigma)=1;
X_\sigma (\r)=\emptyset,
\delta (\sigma)=0,
X_{\tau\sigma}(\r)=T_1\amalg 4T_0,
\delta (\tau\sigma )=1,
$}
\centerline{
$
Y(\r)_\sigma=\emptyset,
Y(\r)_{\tau\sigma}=U_1\amalg 2T_0,
Y(\r)=U_1\amalg 2T_0;
$}

\smallpagebreak

\centerline{
$
\alpha (\sigma)=0;
X_\sigma (\r)=\emptyset,
\delta (\sigma)=0,
X_{\tau\sigma}(\r)=T_2\amalg 5T_0,
\delta (\tau\sigma )=1,
$}
\centerline{
$
Y(\r)_\sigma=\emptyset,
Y(\r)_{\tau\sigma}=U_2\amalg U_0\amalg 2T_0,
Y(\r)=U_2\amalg U_0\amalg 2T_0;
$}

\smallpagebreak

\centerline{
$
\alpha (\sigma)=0;
X_\sigma (\r)=2T_1,
\delta (\sigma)=0,
X_{\tau\sigma}(\r)=4T_0,
\delta (\tau\sigma )=1,
$}
\centerline{
$
(Y(\r)_\sigma,  Y(\r)_{\tau\sigma}, Y(\r))=
(2U_1, 2T_0, 2U_1\amalg 2T_0)\ or\
(T_1, 2U_0\amalg T_0, 2U_0\amalg T_1\amalg T_0);
$}

\smallpagebreak

\centerline{
$
\alpha (\sigma)=1;
X_\sigma (\r)=2T_1,
\delta (\sigma)=0,
X_{\tau\sigma}(\r)=4T_0,
\delta (\tau\sigma )=1,
$}
\centerline{
$
Y(\r)_\sigma=T_1,
Y(\r)_{\tau\sigma}=2T_0,
Y(\r)=T_1\amalg 2T_0;
$}

\smallpagebreak

\centerline{
$
\alpha (\sigma)=0;
X_\sigma (\r)=2T_1,
\delta (\sigma)=0,
X_{\tau\sigma}(\r)=T_1\amalg 4T_0,
\delta (\tau\sigma )=1,
$}
\centerline{
$
Y(\r)_\sigma=T_1,
Y(\r)_{\tau\sigma}=U_1\amalg 2T_0,
Y(\r)=U_1\amalg T_1\amalg 2T_0;
$}

\smallpagebreak

\centerline{
$
\alpha (\sigma)=0;
X_\sigma (\r)=T_3\amalg 2T_0,
\delta (\sigma)=0,
X_{\tau\sigma}(\r)=4T_0,
\delta (\tau\sigma )=1,
$}
\centerline{
$
Y(\r)_\sigma=U_3\amalg T_0,
Y(\r)_{\tau\sigma}=2T_0,
Y(\r)=U_3\amalg 3T_0.
$}

\bigpagebreak

The case $(r(\theta), a(\theta ), \delta (\theta ))=(7,1,1)$:

\smallpagebreak

\centerline{
$
\alpha (\sigma)=1;
X_\sigma (\r)=\emptyset,
\delta (\sigma)=0,
X_{\tau\sigma}(\r)=6T_0,
\delta (\tau\sigma )=1,
$}
\centerline{
$
Y(\r)_\sigma=\emptyset,
Y(\r)_{\tau\sigma}=Y(\r)=2U_0\amalg 2T_0\ or\ 3T_0;
$}

\smallpagebreak

\centerline{
$
\alpha (\sigma)=0;
X_\sigma (\r)=\emptyset,
\delta (\sigma)=0,
X_{\tau\sigma}(\r)=6T_0,
\delta (\tau\sigma )=1,
$}
\centerline{
$
Y(\r)_\sigma=\emptyset,
Y(\r)_{\tau\sigma}=Y(\r)=4U_0\amalg T_0\ or\ 2U_0\amalg 2T_0\
or\ 3T_0;
$}

\smallpagebreak

\centerline{
$
\alpha (\sigma)=0;
X_\sigma (\r)=\emptyset,
\delta (\sigma)=0,
X_{\tau\sigma}(\r)=T_1\amalg 6T_0,
\delta (\tau\sigma )=1,
$}
\centerline{
$
Y(\r)_\sigma=\emptyset,
Y(\r)_{\tau\sigma}=Y(\r)=U_1\amalg 2U_0\amalg 2T_0\ or\
U_1\amalg 3T_0;
$}

\smallpagebreak

\centerline{
$
\alpha (\sigma)=1;
X_\sigma (\r)=\emptyset,
\delta (\sigma)=0,
X_{\tau\sigma}(\r)=T_1\amalg 6T_0,
\delta (\tau\sigma )=1,
$}
\centerline{
$
Y(\r)_\sigma=\emptyset,
Y(\r)_{\tau\sigma}=U_1\amalg 3T_0,
Y(\r)=U_1\amalg 3T_0;
$}

\smallpagebreak

\centerline{
$
\alpha (\sigma)=0;
X_\sigma (\r)=\emptyset,
\delta (\sigma)=0,
X_{\tau\sigma}(\r)=T_2\amalg 7T_0,
\delta (\tau\sigma )=1,
$}
\centerline{
$
Y(\r)_\sigma=\emptyset,
Y(\r)_{\tau\sigma}=U_2\amalg U_0\amalg 3T_0,
Y(\r)=U_2\amalg U_0\amalg 3T_0;
$}

\smallpagebreak

\centerline{
$
\alpha (\sigma)=0;
X_\sigma (\r)=2T_1,
\delta (\sigma)=0,
X_{\tau\sigma}(\r)=6T_0,
\delta (\tau\sigma )=1,
$}
\centerline{
$
(Y(\r)_\sigma, Y(\r)_{\tau\sigma}, Y(\r))=(2U_1, 3T_0, 2U_1\amalg 3T_0)\
or\
(T_1, 2U_0\amalg 2T_0, 2U_0\amalg T_1\amalg 2T_0);
$}

\smallpagebreak

\centerline{
$
\alpha (\sigma)=1;
X_\sigma (\r)=2T_1,
\delta (\sigma)=0,
X_{\tau\sigma}(\r)=6T_0,
\delta (\tau\sigma )=1,
$}
\centerline{
$
Y(\r)_\sigma=T_1,
Y(\r)_{\tau\sigma}=3T_0,
Y(\r)=T_1\amalg 3T_0;
$}

\smallpagebreak

\centerline{
$
\alpha (\sigma)=0;
X_\sigma (\r)=2T_1,
\delta (\sigma)=0,
X_{\tau\sigma}(\r)=T_1\amalg 6T_0,
\delta (\tau\sigma )=1,
$}
\centerline{
$
Y(\r)_\sigma=T_1,
Y(\r)_{\tau\sigma}=U_1\amalg 3T_0,
Y(\r)=U_1\amalg T_1\amalg 3T_0;
$}

\smallpagebreak

\centerline{
$
\alpha (\sigma)=0;
X_\sigma (\r)=T_3\amalg 2T_0,
\delta (\sigma)=0,
X_{\tau\sigma}(\r)=6T_0,
\delta (\tau\sigma )=1,
$}
\centerline{
$
Y(\r)_\sigma=U_3\amalg T_0,
Y(\r)_{\tau\sigma}=3T_0,
Y(\r)=U_3\amalg 4T_0;
$}

\smallpagebreak

\centerline{
$
\alpha (\sigma)=0;
X_\sigma (\r)=8T_0,
\delta (\sigma)=0,
X_{\tau\sigma}(\r)=T_3,
\delta (\tau\sigma )=1,
$}
\centerline{
$
Y(\r)_\sigma=4T_0,
Y(\r)_{\tau\sigma}=U_3,
Y(\r)=U_3\amalg 4T_0.
$}

\bigpagebreak

The case $(r(\theta), a(\theta ), \delta (\theta ))=(4,2,0)$:

\smallpagebreak

\centerline{
$
\alpha (\sigma)=1;
X_\sigma (\r)=2T_1,
\delta (\sigma)=0,
X_{\tau\sigma}(\r)=\emptyset,
\delta (\tau\sigma )=0,
$}
\centerline{
$
Y(\r)_{\tau\sigma}=\emptyset,
Y(\r)_\sigma = Y(\r)=2U_1\ or\ T_1,
$}

\smallpagebreak

\centerline{
$
\alpha (\sigma)=0;
X_\sigma (\r)=2T_1,
\delta (\sigma)=0,
X_{\tau\sigma}(\r)=\emptyset,
\delta (\tau\sigma )=0,
$}
\centerline{
$
Y(\r)_{\tau\sigma}=\emptyset,
Y(\r)_\sigma = Y(\r)=2U_1\ or\ T_1;
$}

\smallpagebreak

\centerline{
$
\alpha (\sigma)=0;
X_\sigma (\r)=T_3\amalg 2T_0,
\delta (\sigma)=0,
X_{\tau\sigma}(\r)=\emptyset,
\delta (\tau\sigma )=0,
$}
\centerline{
$
Y(\r)_{\tau\sigma}=\emptyset,
Y(\r)_\sigma = Y(\r)=U_3\amalg 2U_0\ or\ U_3\amalg T_0;
$}

\smallpagebreak

\centerline{
$
\alpha (\sigma)=1;
X_\sigma (\r)=\emptyset,
\delta (\sigma)=0,
X_{\tau\sigma}(\r)=\emptyset,
\delta (\tau\sigma )=0,
$}
\centerline{
$
Y(\r)_\sigma=\emptyset,
Y(\r)_{\tau\sigma}=\emptyset,
Y(\r)=\emptyset;
$}

\smallpagebreak

\centerline{
$
\alpha (\sigma)=1;
X_\sigma (\r)=T_3\amalg 2T_0,
\delta (\sigma)=0,
X_{\tau\sigma}(\r)=\emptyset,
\delta (\tau\sigma )=0,
$}
\centerline{
$
Y(\r)_\sigma=U_3\amalg T_0,
Y(\r)_{\tau\sigma}=\emptyset,
Y(\r)=U_3\amalg T_0;
$}

\smallpagebreak

\centerline{
$
\alpha (\sigma)=0;
X_\sigma (\r)=T_4\amalg 3T_0,
\delta (\sigma)=0,
X_{\tau\sigma}(\r)=\emptyset,
\delta (\tau\sigma )=0,
$}
\centerline{
$
Y(\r)_\sigma=U_4\amalg U_0\amalg T_0,
Y(\r)_{\tau\sigma}=\emptyset,
Y(\r)=U_4\amalg U_0\amalg T_0;
$}

\smallpagebreak

\centerline{
$
\alpha (\sigma)=0;
X_\sigma (\r)=2T_1,
\delta (\sigma)=0,
X_{\tau\sigma}(\r)=2T_1,
\delta (\tau\sigma )=0,
$}
\centerline{
$
(Y(\r)_\sigma , Y(\r)_{\tau\sigma})=(2U_1,T_1)\ or\ (T_1, 2U_1),
Y(\r)=2U_1\amalg T_1;
$}

\smallpagebreak

\centerline{
$
\alpha (\sigma)=1;
X_\sigma (\r)=2T_1,
\delta (\sigma)=0,
X_{\tau\sigma}(\r)=2T_1,
\delta (\tau\sigma )=0,
$}
\centerline{
$
Y(\r)_\sigma=T_1,
Y(\r)_{\tau\sigma}=T_1,
Y(\r)=2T_1;
$}

\smallpagebreak

\centerline{
$
\alpha (\sigma)=0;
X_\sigma (\r)=T_3\amalg 2T_0,
\delta (\sigma)=0,
X_{\tau\sigma}(\r)=2T_1,
\delta (\tau\sigma )=0,
$}
\centerline{
$
Y(\r)_\sigma=U_3\amalg T_0,
Y(\r)_{\tau\sigma}=T_1,
Y(\r)=U_3\amalg T_1\amalg T_0;
$}

\smallpagebreak

\centerline{
$
\alpha (\sigma)=1;
X_\sigma (\r)=T_2,
\delta _{\sigma S}=1,
X_{\tau\sigma}(\r)=T_2\amalg 2T_0,
\delta _{\tau\sigma S}=1,
$}
\centerline{
$
Y(\r)_\sigma=U_2,
Y(\r)_{\tau\sigma}=U_2\amalg T_0,
Y(\r)=2U_2\amalg T_0;
$}

\smallpagebreak

\centerline{
$
\alpha (\sigma)=1;
X_\sigma (\r)=T_4\amalg 2T_0,
\delta _{\sigma S}=1,
X_{\tau\sigma}(\r)=T_0,
\delta _{\tau\sigma S}=1,
$}
\centerline{
$
Y(\r)_\sigma=U_4\amalg T_0,
Y(\r)_{\tau\sigma}=U_0,
Y(\r)=U_4\amalg U_0\amalg T_0;
$}

\smallpagebreak

\centerline{
$
\alpha (\sigma)=1;
X_\sigma (\r)=T_3\amalg T_0,
\delta _{\sigma S}=1,
X_{\tau\sigma}(\r)=T_0,
\delta _{\tau\sigma S}=1,
$}
\centerline{
$
Y(\r)_\sigma=U_3\amalg U_0,
Y(\r)_{\tau\sigma}=U_0,
Y(\r)=U_3\amalg 2U_0;
$}

\smallpagebreak

\centerline{
$
\alpha (\sigma)=1;
X_\sigma (\r)=T_2,
\delta _{\sigma S}=1,
X_{\tau\sigma}(\r)=T_1\amalg T_0,
\delta _{\tau\sigma S}=1,
$}
\centerline{
$
Y(\r)_\sigma=U_2,
Y(\r)_{\tau\sigma}=U_1\amalg U_0,
Y(\r)=U_2\amalg U_1\amalg U_0;
$}

\smallpagebreak

\centerline{
$
\alpha (\sigma)=1;
X_\sigma (\r)=T_2,
\delta _{\sigma S}=1,
X_{\tau\sigma}(\r)=T_2\amalg 2T_0,
\delta _{\tau\sigma S}=1,
$}
\centerline{
$
Y(\r)_\sigma=U_2,
Y(\r)_{\tau\sigma}=U_2\amalg 2U_0,
Y(\r)=2U_2\amalg 2U_0;
$}

\smallpagebreak

\centerline{
$
\alpha (\sigma)=1;
X_\sigma (\r)=T_3\amalg 2T_0,
\delta _{\sigma S}=0,
\delta _{\sigma S_+\cap S_-}=1,
$}
\centerline{
$
X_{\tau\sigma}(\r)=T_1,
\delta _{\tau\sigma S}=0,
\delta _{\tau\sigma S_+\cap S_-}=1,
$}
\centerline{
$
Y(\r)_\sigma=U_3\amalg T_0,
Y(\r)_{\tau\sigma}=U_1,
Y(\r)=U_3\amalg U_1\amalg T_0;
$}

\smallpagebreak

\centerline{
$
\alpha (\sigma)=1;
X_\sigma (\r)=T_3,
\delta _{\sigma S}=0,
\delta _{\sigma S_+\cap S_-}=1,
$}
\centerline{
$
X_{\tau\sigma}(\r)=T_1\amalg 2T_0,
\delta _{\tau\sigma S}=0,
\delta _{\tau\sigma S_+\cap S_-}=1,
$}
\centerline{
$
Y(\r)_\sigma=U_3,
Y(\r)_{\tau\sigma}=U_1\amalg T_0,
Y(\r)=U_3\amalg U_1\amalg T_0;
$}

\smallpagebreak

\centerline{
$
\alpha (\sigma)=1;
X_\sigma (\r)=T_4\amalg T_0,
\delta _{\sigma S}=0,
\delta _{\sigma S_+\cap S_-}=1,
$}
\centerline{
$
X_{\tau\sigma}(\r)=2T_0,
\delta _{\tau\sigma S}=0,
\delta _{\tau\sigma S_+\cap S_-}=1,
$}
\centerline{
$
Y(\r)_\sigma=U_4\amalg U_0,
Y(\r)_{\tau\sigma}=T_0,
Y(\r)=U_4\amalg U_0\amalg T_0;
$}

\smallpagebreak

\centerline{
$
\alpha (\sigma)=1;
X_\sigma (\r)=T_2\amalg T_0,
\delta _{\sigma S}=0,
\delta _{\sigma S_+\cap S_-}=1,
$}
\centerline{
$
X_{\tau\sigma}(\r)=T_1,
\delta _{\tau\sigma S}=0,
\delta _{\tau\sigma S_+\cap S_-}=1,
$}
\centerline{
$
Y(\r)_\sigma=U_2\amalg U_0,
Y(\r)_{\tau\sigma}=U_1,
Y(\r)=U_2\amalg U_1\amalg U_0;
$}

\smallpagebreak

\centerline{
$
\alpha (\sigma)=1;
X_\sigma (\r)=T_3,
\delta _{\sigma S}=0,
\delta _{\sigma S_+\cap S_-}=1,
$}
\centerline{
$
X_{\tau\sigma}(\r)=2T_0,
\delta _{\tau\sigma S}=0,
\delta _{\tau\sigma S_+\cap S_-}=1,
$}
\centerline{
$
Y(\r)_\sigma=U_3,
Y(\r)_{\tau\sigma}=2U_0,
Y(\r)=U_3\amalg 2U_0;
$}

\smallpagebreak

\centerline{
$
\alpha (\sigma)=0;
X_\sigma (\r)=T_3,
\delta (\sigma )=1,
\delta _{\sigma S_+\cap S_-}=0,
$}
\centerline{
$
X_{\tau\sigma}(\r)=4T_0,
\delta (\tau\sigma )=1,
\delta _{\tau\sigma S_+\cap S_-}=0,
$}
\centerline{
$
Y(\r)_\sigma = U_3,
Y(\r)_{\tau\sigma} = 2T_0,
Y(\r)=U_3\amalg 2T_0.
$}

\bigpagebreak

The case $(r(\theta), a(\theta ), \delta (\theta ))=(8,0,0)$:

\smallpagebreak

\centerline{
$
\alpha (\sigma)=1;
X_\sigma (\r)=8T_0,
\delta (\sigma)=0,
X_{\tau\sigma}(\r)=\emptyset,
\delta (\tau\sigma )=0,
$}
\centerline{
$
Y(\r)_{\tau\sigma}=\emptyset,
Y(\r)_\sigma = Y(\r)=2U_0\amalg 3T_0\ or\ 4T_0;
$}

\smallpagebreak

\centerline{
$
\alpha (\sigma)=0;
X_\sigma (\r)=8T_0,
\delta (\sigma)=0,
X_{\tau\sigma}(\r)=\emptyset,
\delta (\tau\sigma )=0,
$}
\centerline{
$
Y(\r)_{\tau\sigma}=\emptyset,
Y_\sigma (\r)=Y(\r)=4U_0\amalg 2T_0\ or\ 2U_0\amalg 3T_0\ or\ 4T_0;
$}

\smallpagebreak

\centerline{
$
\alpha (\sigma)=0;
X_\sigma (\r)=T_1\amalg 8T_0,
\delta (\sigma)=0,
X_{\tau\sigma}(\r)=\emptyset,
\delta (\tau\sigma )=0,
$}
\centerline{
$
Y(\r)_{\tau\sigma}=\emptyset,
Y(\r)_\sigma=Y(\r)=U_1\amalg 2U_0\amalg 3T_0\ or\ U_1\amalg 4T_0;
$}

\smallpagebreak

\centerline{
$
\alpha (\sigma)=1;
X_\sigma (\r)=T_1\amalg 8T_0,
\delta (\sigma)=0,
X_{\tau\sigma}(\r)=\emptyset,
\delta (\tau\sigma )=0,
$}
\centerline{
$
Y(\r)_\sigma=U_1\amalg 4T_0,
Y(\r)_{\tau\sigma}=\emptyset,
Y(\r)=U_1\amalg 4T_0;
$}

\smallpagebreak

\centerline{
$
\alpha (\sigma)=0;
X_\sigma (\r)=T_2\amalg 9T_0,
\delta (\sigma)=0,
X_{\tau\sigma}(\r)=\emptyset,
\delta (\tau\sigma )=0,
$}
\centerline{
$
Y(\r)_\sigma=U_2\amalg U_0\amalg 4T_0,
Y(\r)_{\tau\sigma}=\emptyset,
Y(\r)=U_2\amalg U_0\amalg 4T_0;
$}

\smallpagebreak

\centerline{
$
\alpha (\sigma)=0;
X_\sigma (\r)=8T_0,
\delta (\sigma)=0,
X_{\tau\sigma}(\r)=2T_1,
\delta (\tau\sigma )=0,
$}
\centerline{
$
(Y(\r)_\sigma, Y(\r)_{\tau\sigma}, Y(\r))=
(4T_0, 2U_1, 2U_1\amalg 4T_0)\ or\
(2U_0\amalg 3T_0, T_1, 2U_0\amalg T_1 \amalg 3T_0);
$}

\smallpagebreak

\centerline{
$
\alpha (\sigma)=1;
X_\sigma (\r)=8T_0,
\delta (\sigma)=0,
X_{\tau\sigma}(\r)=2T_1,
\delta (\tau\sigma )=0,
$}
\centerline{
$
Y(\r)_\sigma=4T_0,
Y(\r)_{\tau\sigma}=T_1,
Y(\r)=T_1\amalg 4T_0;
$}

\smallpagebreak

\centerline{
$
\alpha (\sigma)=0;
X_\sigma (\r)=T_1\amalg 8T_0,
\delta (\sigma)=0,
X_{\tau\sigma}(\r)=2T_1,
\delta (\tau\sigma )=0,
$}
\centerline{
$
Y(\r)_\sigma=U_1\amalg 4T_0,
Y(\r)_{\tau\sigma}=T_1,
Y(\r)=U_1\amalg T_1\amalg 4T_0;
$}

\smallpagebreak

\centerline{
$
\alpha (\sigma)=0;
X_\sigma (\r)=8T_0,
\delta (\sigma)=0,
X_{\tau\sigma}(\r)=T_3\amalg 2T_0,
\delta (\tau\sigma )=0,
$}
\centerline{
$
Y(\r)_\sigma=4T_0,
Y(\r)_{\tau\sigma}=U_3\amalg T_0,
Y(\r)=U_3\amalg 5T_0;
$}

\smallpagebreak

\centerline{
$
\alpha (\sigma)=1;
X_\sigma (\r)=5T_0,
\delta _{\sigma S}=1,
X_{\tau\sigma}(\r)=3T_0,
\delta _{\tau\sigma S}=1,
$}
\centerline{
$
(Y(\r)_\sigma, Y(\r)_{\tau\sigma})=
(U_0\amalg 2T_0, 3U_0)\ or\ (3U_0\amalg T_0, U_0\amalg T_0),
Y(\r)=4U_0\amalg 2T_0;
$}

\smallpagebreak

\centerline{
$
\alpha (\sigma)=1;
X_\sigma (\r)=5T_0,
\delta _{\sigma S}=1,
X_{\tau\sigma}(\r)=T_2\amalg 4T_0,
\delta _{\tau\sigma S}=1,
$}
\centerline{
$
Y(\r)_\sigma=U_0\amalg 2T_0,
Y(\r)_{\tau\sigma}=U_2\amalg 2T_0,
Y(\r)=U_2\amalg U_0\amalg 4T_0;
$}

\smallpagebreak

\centerline{
$
\alpha (\sigma)=1;
X_\sigma (\r)=T_2\amalg 6T_0,
\delta _{\sigma S}=1,
X_{\tau\sigma}(\r)=3T_0,
\delta _{\tau\sigma S}=1,
$}
\centerline{
$
Y(\r)_\sigma=U_2\amalg 3T_0,
Y(\r)_{\tau\sigma}=U_0\amalg T_0,
Y(\r)=U_2\amalg U_0\amalg 4T_0;
$}

\smallpagebreak

\centerline{
$
\alpha (\sigma)=1;
X_\sigma (\r)=7T_0,
\delta _{\sigma S}=1,
X_{\tau\sigma}(\r)=T_2\amalg 2T_0,
\delta _{\tau\sigma S}=1,
$}
\centerline{
$
Y(\r)_\sigma=U_0\amalg 3T_0,
Y(\r)_{\tau\sigma}=U_2\amalg T_0,
Y(\r)=U_2\amalg U_0\amalg 4T_0;
$}

\smallpagebreak

\centerline{
$
\alpha (\sigma)=1;
X_\sigma (\r)=T_1\amalg 5T_0,
\delta _{\sigma S}=1,
X_{\tau\sigma}(\r)=3T_0,
\delta _{\tau\sigma S}=1,
$}
\centerline{
$
Y(\r)_\sigma=U_1\amalg U_0\amalg 2T_0,
Y(\r)_{\tau\sigma}=U_0\amalg T_0,
Y(\r)=U_1\amalg 2U_0\amalg 3T_0;
$}

\smallpagebreak

\centerline{
$
\alpha (\sigma)=1;
X_\sigma (\r)=5T_0,
\delta _{\sigma S}=1,
X_{\tau\sigma}(\r)=T_1\amalg 3T_0,
\delta _{\tau\sigma S}=1,
$}
\centerline{
$
Y(\r)_\sigma=U_0\amalg 2T_0,
Y(\r)_{\tau\sigma}=U_1\amalg U_0\amalg T_0,
Y(\r)=U_1\amalg 2U_0\amalg 3T_0;
$}

\smallpagebreak

\centerline{
$
\alpha (\sigma)=1;
X_\sigma (\r)=T_1\amalg 7T_0,
\delta _{\sigma S}=1,
X_{\tau\sigma}(\r)=T_0,
\delta _{\tau\sigma S}=1,
$}
\centerline{
$
Y(\r)_\sigma=U_1\amalg U_0\amalg 3T_0,
Y(\r)_{\tau\sigma}=U_0,
Y(\r)=U_1\amalg 2U_0\amalg 3T_0;
$}

\smallpagebreak

\centerline{
$
\alpha (\sigma)=1;
X_\sigma (\r)=7T_0,
\delta _{\sigma S}=1,
X_{\tau\sigma}(\r)=T_1\amalg T_0,
\delta_{\tau\sigma S}=1,
$}
\centerline{
$
Y(\r )_\sigma =U_0\amalg 3T_0,
Y(\r )_{\tau\sigma}=U_1\amalg U_0,
Y(\r )=U_1\amalg 2U_0\amalg 3T_0;
$}

\smallpagebreak

\centerline{
$
\alpha (\sigma)=1;
X_\sigma (\r)=7T_0,
\delta _{\sigma S}=1,
X_{\tau\sigma}(\r)=T_0,
\delta_{\tau\sigma S}=1,
$}
\centerline{
$
Y(\r )_\sigma =3U_0\amalg 2T_0,
Y(\r )_{\tau\sigma}=U_0,
Y(\r )=4U_0\amalg 2T_0;
$}

\smallpagebreak

\centerline{
$
\alpha (\sigma)=1;
X_\sigma (\r)=T_1\amalg 4T_0,
\delta _{\sigma S}=0,
\delta_{\sigma S_+\cap S_-}=1,
$}
\centerline{
$
X_{\tau\sigma}(\r)=4T_0,
\delta_{\tau\sigma S}=0,
\delta_{\tau\sigma S_+\cap S_-}=1,
$}
\centerline{
$
(Y(\r )_\sigma , Y(\r )_{\tau\sigma})=
(U_1\amalg 2U_0\amalg T_0,2T_0)\ or\ (U_1\amalg 2T_0, 2U_0\amalg T_0),
$}
\centerline{
$
Y(\r )=U_1\amalg 2U_0\amalg 3T_0;
$}

\smallpagebreak

\centerline{
$
\alpha (\sigma)=1;
X_\sigma (\r)=T_1\amalg 6T_0,
\delta _{\sigma S}=0,
\delta_{\sigma S_+\cap S_-}=1,
$}
\centerline{
$
X_{\tau\sigma}(\r)=2T_0,
\delta_{\tau\sigma S}=0,
\delta_{\tau\sigma S_+\cap S_-}=1,
$}
\centerline{
$
(Y(\r )_\sigma , Y(\r )_{\tau\sigma})=
(U_1\amalg 2U_0\amalg 2T_0 ,T_0)\ or\ (U_1\amalg 3T_0,2U_0),
$}
\centerline{
$
Y(\r )=U_1\amalg 2U_0\amalg 3T_0;
$}

\smallpagebreak

\centerline{
$
\alpha (\sigma)=1;
X_\sigma (\r)=6T_0,
\delta _{\sigma S}=0,
\delta_{\sigma S_+\cap S_-}=1,
$}
\centerline{
$
X_{\tau\sigma}(\r)=T_1\amalg 2T_0,
\delta_{\tau\sigma S}=0,
\delta_{\tau\sigma S_+\cap S_-}=1,
$}
\centerline{
$
(Y(\r )_\sigma , Y(\r )_{\tau\sigma})=
(2U_0\amalg 2T_0, U_1\amalg T_0)\ or\ (3T_0, U_1\amalg 2U_0),
$}
\centerline{
$
Y(\r )=U_1\amalg 2U_0\amalg 3T_0;
$}

\smallpagebreak

\centerline{
$
\alpha (\sigma)=1;
X_\sigma (\r)=4T_0,
\delta _{\sigma S}=0,
\delta_{\sigma S_+\cap S_-}=1,
$}
\centerline{
$
X_{\tau\sigma}(\r)=4T_0,
\delta_{\tau\sigma S}=0,
\delta_{\tau\sigma S_+\cap S_-}=1,
$}
\centerline{
$
(Y(\r )_\sigma , Y(\r )_{\tau\sigma})=
(2U_0\amalg T_0, 2U_0\amalg T_0)\ or\ (4U_0,2T_0),
$}
\centerline{
$
Y(\r )=4U_0\amalg 2T_0;
$}

\smallpagebreak

\centerline{
$
\alpha (\sigma)=1;
X_\sigma (\r)=6T_0,
\delta _{\sigma S}=0,
\delta_{\sigma S_+\cap S_-}=1,
$}
\centerline{
$
X_{\tau\sigma}(\r)=2T_0,
\delta_{\tau\sigma S}=0,
\delta_{\tau\sigma S_+\cap S_-}=1,
$}
\centerline{
$
(Y(\r )_\sigma , Y(\r )_{\tau\sigma})=
(4U_0\amalg T_0,T_0)\ or\ (2U_0\amalg 2T_0,2U_0),
$}
\centerline{
$
Y(\r )=4U_0\amalg 2T_0;
$}

\smallpagebreak

\centerline{
$
\alpha (\sigma )=1; X_\sigma (\r )=T_1\amalg 4T_0,
\delta_{\sigma S}=0,
\delta _{\sigma S_+\cap S_-}=1,
$}
\centerline{
$
X_{\tau\sigma}(\r)=T_1\amalg 4T_0,
\delta _{\tau\sigma S}=0,
\delta _{\tau\sigma S_+\cap S_-}=1,
$}
\centerline{
$
Y(\r)_\sigma=U_1\amalg 2T_0,
Y(\r)_{\tau\sigma}=U_1\amalg 2T_0,
Y(\r)=2U_1\amalg 4T_0;
$}

\smallpagebreak

\centerline{
$
\alpha (\sigma)=1;
X_\sigma (\r)=T_1\amalg 6T_0,
\delta _{\sigma S}=0,
\delta _{\sigma S_+\cap S_-}=1,
$}
\centerline{
$
X_{\tau\sigma}(\r)=T_1\amalg 2T_0,
\delta _{\tau\sigma S}=0,
\delta _{\tau\sigma S_+\cap S_-}=1,
$}
\centerline{
$
Y(\r)_\sigma=U_1\amalg 3T_0,
Y(\r)_{\tau\sigma}=U_1\amalg T_0,
Y(\r)=2U_1\amalg 4T_0;
$}

\smallpagebreak

\centerline{
$
\alpha (\sigma)=1;
X_\sigma (\r)=T_1\amalg 8T_0,
\delta _{\sigma S}=0,
\delta _{\sigma S_+\cap S_-}=1,
$}
\centerline{
$
X_{\tau\sigma}(\r)=T_1,
\delta _{\tau\sigma S}=0,
\delta _{\tau\sigma S_+\cap S_-}=1,
$}
\centerline{
$
Y(\r)_\sigma =U_1\amalg 4T_0,
Y(\r)_{\tau\sigma}=U_1,
Y(\r)=2U_1\amalg 4T_0;
$}

\smallpagebreak

\centerline{
$
\alpha (\sigma)=1;
X_\sigma (\r)=T_2\amalg 5T_0,
\delta _{\sigma S}=0,
\delta _{\sigma S_+\cap S_-}=1,
$}
\centerline{
$
X_{\tau\sigma}(\r)=4T_0,
\delta _{\tau\sigma S}=0,
\delta _{\tau\sigma S_+\cap S_-}=1,
$}
\centerline{
$
Y(\r)_\sigma=U_2\amalg U_0\amalg 2T_0,
Y(\r)_{\tau\sigma}=2T_0,
Y(\r)=U_2\amalg U_0\amalg 4T_0;
$}

\smallpagebreak

\centerline{
$
\alpha (\sigma)=1;
X_\sigma (\r)=T_2\amalg 7T_0,
\delta _{\sigma S}=0,
\delta _{\sigma S_+\cap S_-}=1,
$}
\centerline{
$
X_{\tau\sigma}(\r)=2T_0,
\delta _{\tau\sigma S}=0,
\delta _{\tau\sigma S_+\cap S_-}=1,
$}
\centerline{
$
Y(\r)_\sigma=U_2\amalg U_0\amalg 3T_0,
Y(\r)_{\tau\sigma}=T_0,
Y(\r)=U_2\amalg U_0\amalg 4T_0;
$}

\smallpagebreak

\centerline{
$
\alpha (\sigma)=1;
X_\sigma (\r)=6T_0,
\delta _{\sigma S}=0,
\delta _{\sigma S_+\cap S_-}=1,
$}
\centerline{
$
X_{\tau\sigma}(\r)=T_2\amalg 3T_0,
\delta _{\tau\sigma S}=0,
\delta _{\tau\sigma S_+\cap S_-}=1,
$}
\centerline{
$
Y(\r)_\sigma=3T_0,
Y(\r)_{\tau\sigma}=U_2\amalg U_0\amalg T_0,
Y(\r)=U_2\amalg U_0\amalg 4T_0;
$}

\smallpagebreak

\centerline{
$
\alpha (\sigma)=1;
X_\sigma (\r)=8T_0,
\delta _{\sigma S}=0,
\delta _{\sigma S_+\cap S_-}=1,
$}
\centerline{
$
X_{\tau\sigma}(\r)=T_2\amalg T_0,
\delta _{\tau\sigma S}=0,
\delta _{\tau\sigma S_+\cap S_-}=1,
$}
\centerline{
$
Y(\r)_\sigma=4T_0,
Y(\r)_{\tau\sigma}=U_2\amalg U_0,
Y(\r)=U_2\amalg U_0\amalg 4T_0;
$}

\smallpagebreak

\centerline{
$
\alpha (\sigma)=1;
X_\sigma (\r)=8T_0,
\delta _{\sigma S}=0,
\delta _{\sigma S_+\cap S_-}=1,
$}
\centerline{
$
X_{\tau\sigma}(\r)=T_1,
\delta _{\tau\sigma S}=0,
\delta _{\tau\sigma S_+\cap S_-}=1,
$}
\centerline{
$
Y(\r)_\sigma=2U_0\amalg 3T_0,
Y(\r)_{\tau\sigma}=U_1,
Y(\r)=U_1\amalg 2U_0\amalg 3T_0.
$}

\bigpagebreak

As a result, we get the following

\proclaim{Theorem 6} There exist the following $59$ topological
types of $Y(\r)$ for real Enriques surfaces $Y$:

$s_{nor}=0:
\ \ \emptyset;\ \ \ T_1\amalg T_0,\ 2T_1;\ \ \ T_1\amalg 2T_0;\ \ \
4T_0,\ T_1\amalg 3T_0;\ \ \ 5T_0,\ T_1\amalg 4T_0.$

$s_{nor}=1:\ \
U_{1+2k}, k=0,1,2,3,4;\ \ \
U_1\amalg T_0,\ U_1\amalg T_1,\
U_3\amalg T_0,\ U_3\amalg T_1,\ U_5\amalg T_1,\ U_7\amalg T_1,\
U_9\amalg T_0,\ U_9\amalg T_1;\ \ \
U_1\amalg 2T_0,\ U_1\amalg T_1\amalg T_0,\ U_3\amalg 2T_0,\
U_3\amalg T_1\amalg T_0;\ \ \
U_1\amalg 3T_0,\ U_1\amalg T_1\amalg 2T_0,\ U_3\amalg 3T_0;\ \ \
U_1\amalg 4T_0,\ U_1\amalg T_1\amalg 3T_0,\ U_3\amalg 4T_0;\ \ \
U_1\amalg 5T_0,\ U_1\amalg T_1\amalg 4T_0;\ \ \
U_3\amalg 5T_0. $

$s_{nor}=2:\ \
U_{2k}\amalg U_0, k=1,2,3,4,5,\ \ U_{10-k}\amalg U_k,\ k=1,2,3,4,5;\ \ \
2U_1\amalg T_1,\ U_2\amalg U_0\amalg T_0,\ 2U_2\amalg T_0,\
U_3\amalg U_1\amalg T_0,\ U_4\amalg U_0\amalg T_0;\ \ \
U_2\amalg U_0\amalg 2T_0; \ \ \
2U_0\amalg 3T_0,\ U_2\amalg U_0\amalg 3T_0;\ \ \
2U_1\amalg 4T_0,\ U_2\amalg U_0\amalg 4T_0.$

$s_{nor}=3:\ \
U_2\amalg U_1\amalg U_0,\ U_3\amalg 2U_0;\ \ \
U_1\amalg 2U_0\amalg 3T_0. $

$s_{nor}=4:\ \
2U_2\amalg 2U_0;\ \ \ 4U_0\amalg 2T_0. $

\endproclaim

\bigpagebreak

\proclaim{Theorem 7} All other topological types of $Y(\r)$ for
real Enriques surfaces $Y$, different from the types of Theorem 6
belong to the following list of $21$ topological types:

$s_{nor}=0:\ \ T_0,\ T_1;\ \ \ 2T_0;\ \ \ 3T_0. $

$s_{nor}=2:\ \ 2U_0, 2U_1;\ \ \ 2U_0\amalg T_0,\
2U_0\amalg T_1,\ 2U_1\amalg T_0;\ \ \
2U_0\amalg 2T_0,\ 2U_0\amalg T_1\amalg T_0,\ 2U_1\amalg 2T_0;\ \ \
2U_0\amalg T_1\amalg 2T_0,\ 2U_1\amalg 3T_0;\ \ \
2U_0\amalg 4T_0,\ 2U_0\amalg T_1\amalg 3T_0. $

$s_{nor}=3:\ \ U_1\amalg 2U_0;\ \ \ U_1\amalg 2U_0\amalg T_0;\ \ \
U_1\amalg 2U_0\amalg 2T_0. $

$s_{nor}=4:\ \ 4U_0;\ \ \ 4U_0\amalg T_0. $
\endproclaim

\bigpagebreak

We think that calculations above (Theorem 5) of additional invariants
 will be useful to prohibit some
topological types of Theorem 7.

It is very interesting to calculate arithmetic invariants of real
Enriques surfaces, for example the Brauer group and the
Witt group, and compare these invariants with similar
invariants for real
rational surfaces. See \cite{CT--P}, \cite{Kr},
\cite{Ma}, \cite{Mi}, \cite{N5}, \cite{N--S}, \cite{Si}, \cite{Su}
on this subject. We hope that our results here will be
useful for calculations of these invariants.

\enddocument

\bye